\documentclass[11pt,twoside,english]{amsart}
\advance\oddsidemargin by -1.0cm
\advance\evensidemargin by -1.0cm
\advance\topmargin by -1.0cm
\usepackage{amssymb}
\usepackage{babel}
\usepackage{amstext}
\usepackage{amscd}   
\usepackage{epsfig}  
\usepackage{rotating}

\theoremstyle{plain}
\newtheorem{cor}{Corollary}[section]
\newtheorem{lem}{Lemma}[section]
\newtheorem{thm}{Theorem}[section]            

\theoremstyle{definition}

\newtheorem{NB}{Remark}[section]

\newcommand{\bdm}{\begin{displaymath}}
\newcommand{\edm}{\end{displaymath}}
\newcommand{\be}{\begin{equation}}
\newcommand{\ee}{\end{equation}}
\newcommand{\ba}[1]{\begin{array}{#1}}
\newcommand{\ea}{\end{array}}

\newcommand{\btab}{\begin{tabular}}
\newcommand{\etab}{\end{tabular}}

\begin{document}
\def\haken{\mathbin{\hbox to 6pt{%
                 \vrule height0.4pt width5pt depth0pt
                 \kern-.4pt
                 \vrule height6pt width0.4pt depth0pt\hss}}}
    \let \hook\intprod
\setcounter{equation}{0}
%
%
\thispagestyle{empty}
%
\date{\today}
\title[Eigenvalues estimates in terms of Codazzi tensors]
{Eigenvalues estimates for the Dirac operator in terms of  Codazzi tensors}
%
%
%
\author{Th. Friedrich}
\author{E.C. Kim}
\address{\hspace{-5mm} 
{\normalfont\ttfamily friedric@mathematik.hu-berlin.de}\newline
Institut f\"ur Mathematik \newline
Humboldt-Universit\"at zu Berlin\newline
Sitz: WBC Adlershof\newline
D-10099 Berlin, Germany}
\address{\hspace{-5mm} 
{\normalfont\ttfamily eckim@andong.ac.kr}\newline
Andong National University\newline
Department of Mathematics\newline
College of Education\newline
388 Songchon-dong, Andong\newline
760-749 Kyoungsangbuk-do\newline
South Korea}
\thanks{Supported by the SFB 647 "Raum, Zeit,Materie" of the DFG}
\subjclass[2000]{Primary 53 C 25; Secondary 53 C 27}
%
\begin{abstract}
We prove a lower bound for the first eigenvalue of the Dirac operator on a 
compact Riemannian spin manifold depending on the scalar curvature as well as 
a chosen Codazzi tensor. The inequality generalizes the classical
estimate from \cite{Friedrich80}.
%
\end{abstract}
\maketitle
\pagestyle{headings}
%
%
\section{Introduction}
%
The first author proved in \cite{Friedrich80} that the smallest eigenvalue 
$\lambda_1$  of the Dirac operator $D$ of a compact Riemannian spin manifold  
$(M^n, g)$ satisfies
\begin{equation}
\lambda_1^2    \, \geq  \, \frac{n}{4(n-1)} \cdot S_{\mathrm{min}} \ ,
\end{equation}
where $S_{\mathrm{min}}$ denotes the minimum of the scalar curvature.
The limiting case of (1) occurs if and only if $(M^n, g)$ admits a
nontrivial spinor field $\psi_1$ satisfying
\[  \nabla_X \psi_1 =  - \frac{\lambda_1}{n}  \, X  \cdot  \psi_1 , \]
where $X$ is an arbitrary vector field on $M^n$ and the dot "$\cdot $"
indicates Clifford multiplication \cite{Friedrich2000}. Improvements of 
this estimate do typically depend on additional geometric structures on
the considered manifold $(M^n, g)$ 
\cite{FriedrichKim}, \cite{Kirchberg}, \cite{KSW}. 
The aim of this paper is to show that  inequality (1) can be
improved in case that a Codazzi tensor exists.

\medskip
A symmetric $(0, 2)$-tensor field $\beta$ on $(M^n, g )$ is called  a
\emph{nondegenerate Codazzi tensor} \cite{Besse} if  $\beta$ is nondegenerate 
at all points of $M^n$ and satisfies
$  ( \nabla_X  \beta ) (Y, Z) = ( \nabla_Y  \beta ) (X, Z) $  for all vector
fields $X, Y, Z$. We 
identify $\beta$ with the induced $(1,1)$-tensor $\beta$ via 
$\beta(X, Y) = g(X, \beta(Y))$.                                                                     
Let $(E_1, \ldots, E_n )$ be  a local orthonormal frame field on $(M^n, g)$. 
Then the spin derivative $\nabla$ and the Dirac operator $D$, acting on
sections $\psi \in \Gamma( \Sigma(M^n))$ of the spinor bundle $\Sigma(M^n)$
over $(M^n, g)$, are locally expressed as \cite{Friedrich2000}
\begin{eqnarray*}
\nabla_X \psi = X(\psi) + \frac{1}{4} \sum_{i=1}^n E_i \cdot \nabla_X E_i
\cdot \psi \ , \quad D \psi = \sum_{i=1}^n E_i \cdot  \nabla_{E_i} \psi  ,
\end{eqnarray*}
respectively. Moreover, we define the $\beta$-{\it twist} $D_{\beta}$ of 
the Dirac operator $D$ by
\[   D_{\beta} \psi  = \sum_{i=1}^n  \beta^{-1} (E_i) \cdot  \nabla_{E_i} \psi 
= \sum_{i=1}^n  E_i \cdot  \nabla_{\beta^{-1}(E_i)}  \psi  .  \]
\begin{thm} 
Let $(M^n, g)$ be an $n$-dimensional closed  Riemannian spin
manifold and consider a nondegenerate Codazzi tensor $\beta$. Denote by
$\overline{g}$ the metric induced 
by $\beta$ via $\overline{g} (X, Y) = g ( \beta(X), \beta(Y) )$.  
Let $\lambda_1 \in {\mathbb R} $ and $\overline{\lambda}_1 \in {\mathbb R} $
be the smallest eigenvalue  of the Dirac operators $D$ and $\overline{D}$, 
respectively. We assume that both $\lambda_1$ and  $\overline{\lambda}_1$ 
are nonzero.   Then we have 
\begin{equation}
\lambda_1^2 \, \geq \, \inf_M \Big\{ \frac{S}{4(p+1)}  -
\frac{q \, \overline{\lambda}_1^2}{p+1} + \frac{\triangle F}{2 (p+1)F}\Big\},
\end{equation}
where $F : M^n \longrightarrow \mathbb{R}$ is a real-valued
function defined by  
\begin{equation}  
F = - \ \frac{ \vert \det(\beta^{-1} )  \vert }{q} ,  \end{equation} 
$\triangle F := - ( {\rm div}
\circ {\rm grad} )(F)$,   and  $p, q : M^n \longrightarrow
\mathbb{R}$ are bounded real-valued functions satisfying
\begin{equation} 
- \frac{1}{n}  < p < 0 , \quad - \frac{1}{ \vert \beta^{-1} \vert^2} < q < 0,
\end{equation}
that solve the system of two linear equations
\begin{equation}
   n p  +  c ( {\rm tr} \beta^{-1}) q  =  - 1 ,      \qquad
   ( {\rm tr} \beta^{-1}) p +   c \vert \beta^{-1} \vert^2 q =  -  c
\end{equation}
for some nonzero constant $c \not=0  \in \mathbb{R}$.
\end{thm}
The limiting case of (2) occurs  if and only if there exists a spinor field 
$\psi_1$ on $(M^n, g)$ such that 
\begin{equation}
D \psi_1 = \lambda_1 \psi_1 , \qquad D_{\beta} \psi_1  =  \overline{\lambda}_1 \psi_1
\end{equation}
 and
\begin{equation}  
\nabla_X \psi_1 =   \lambda_1 p \, X \cdot \psi_1  +
\overline{\lambda}_1 q \, \beta^{-1}(X) \cdot \psi_1   
\end{equation} 
hold for some nonzero constants$\lambda_1 \not= 0, \overline{\lambda}_1 \not= 0
\in {\mathbb R}$ and for all vector fields $X$. In the limiting case, 
the parameter $c =  
\overline{\lambda}_1/\lambda_1$ is the ratio of the two eigenvalues.\\

If $\beta = g = I$ is the identity map  and $p + q = - 1/n$, then (2) 
reduces to the inequality (1). If the
eigenvalues of $\beta \neq I$ are  constant, but not equal, 
the solutions $p,q$ of the linear system are constant, too,
\begin{eqnarray*}
p(c) \ = \  \frac{\vert \beta^{-1} \vert^2  \, - \,  c \, {\rm tr} \beta^{-1}}{ ( {\rm tr}
  \beta^{-1})^2 \, - \, n \, \vert \beta^{-1} \vert^2} \ , \quad
q(c) \ = \ \frac{ c \, n \, - \,  {\rm tr} \beta^{-1}}{ c ( ( {\rm tr}
  \beta^{-1})^2 \, - \, n \, \vert \beta^{-1} \vert^2)} \ .
\end{eqnarray*}
Consequently, we obtain a family of inequalities depending on a parameter $c
\neq 0$ 
linking $\lambda_1^2 , \, \overline{\lambda}_1^2$ and
$S_{\mathrm{min}}$,
\begin{eqnarray*}
\lambda_1^2  \ + \   \frac{q(c)}{p(c) + 1}\, 
\overline{\lambda}_1^2  \  
\geq   \     \frac{1}{4(p(c) + 1)} \, S_{\mathrm{min}} \ . 
\end{eqnarray*}
The optimal
parameter $c$ is a solution of a quadratic equation, we omit the corresponding
formulas. A 
universal though not optimal value for the parameter $c$ is 
\begin{equation*}
c \ := \ \frac{\vert \beta^{-1} \vert^2}{ {\rm tr} \beta^{-1}}  .
\end{equation*}
In this case we have $ p = 0$ and $q = - 1/\vert \beta^{-1} \vert^2$. This 
particular inequality generalizes (1):
\begin{cor} 
If the eigenvalues of the Codazzi tensor are constant, then
\begin{eqnarray*}
\lambda_1^2  \ \geq   \     \frac{1}{4} \, S_{\mathrm{min}} \ + \ 
 \frac{1}{\vert \beta^{-1} \vert^2 }\, \overline{\lambda}_1^2 \ \geq
\  \frac{1}{4} \, S_{\mathrm{min}} \ + \  
\frac{1}{\vert \beta^{-1} \vert^2} \, 
\frac{n}{4(n-1)} \, \overline{S}_{\mathrm{min}}\, . 
\end{eqnarray*}
\end{cor}
If  ${\rm tr} \beta^{-1} = 0$, the functions $p$ and $q$ do not
depend on the parameter $c$, i.\,e.~we obtain a unique inequality. We will
formulate the result separately. 
\begin{thm} 
Let $(M^n, g)$ be an $n$-dimensional closed  Riemannian spin
manifold and consider a nondegenerate Codazzi tensor such that 
${\rm tr}(\beta^{-1}) = 0$ vanishes identically.
Let $\lambda_1$ and $\overline{\lambda}_1$ be the smallest eigenvalue  of 
$D$ and $\overline{D}$, respectively.  
 We assume that both $\lambda_1$ and  $\overline{\lambda}_1$ are nonzero.  
Then, in the notations of Theorem $1.1$, we have 
\begin{equation}
\lambda_1^2 \, \geq \, \inf_M  \Big\{ \frac{n S}{4(n-1)} +
\frac{n \, \overline{\lambda}_1^2  }{(n-1) \vert \beta^{-1} \vert^2} 
+  \frac{n \, \triangle F}{2 (n-1)F}   \Big\},
\end{equation}
where the real-valued function $F : M^n \longrightarrow \mathbb{R}$ is 
defined by  
\begin{equation}  F =   \vert \det(
\beta^{-1} ) | \cdot  \vert \beta^{-1} \vert^2 .  
\end{equation}        
\end{thm}
The limiting case of (8) occurs  if and only if there exists a spinor field 
$\psi_1$ on $(M^n, g)$ such that 
\begin{equation}  \nabla_X \psi_1 = - \frac{\lambda_1}{n}\, X \cdot \psi_1  -
\frac{\overline{\lambda}_1}{\vert \beta^{-1} \vert^2}  \, \beta^{-1}(X) \cdot
\psi_1   
\end{equation} 
holds
for some constants $\lambda_1 \not= 0 , \overline{\lambda}_1  \not = 0 \in
{\mathbb R}$ and for all vector fields $X$.\\

Let us discuss the $2$-dimensional case in detail. Suppose that $\beta$ is
traceless with eigenvalues $a \, , \, - a$. Then we obtain
\begin{eqnarray*}
 {\rm det}(\beta^{-1}) \ = \ - \, \frac{1}{a^2} \ , \quad  \vert \beta^{-1}
 \vert^2 \ = \ \frac{2}{a^2} \ , \quad F \ \equiv \ \frac{2}{a^4} \ .
\end{eqnarray*}
In particular, the formula of the latter Theorem simplifies:
\begin{cor}
Let $(M^2, g, \beta)$ be a $2$-dimensional closed  Riemannian spin
manifold with a nondegenerate traceless Codazzi tensor. Denote by $\pm a$ its
eigenvalues. Then we have
\begin{eqnarray*}
\lambda_1^2  \,   \geq   \,   \inf_M  \Big\{  \frac{S}{2}   \, + \,
a^2 \,   \overline{\lambda}_1^2 \, + \, a^4 \, \Delta(a^{-4})\Big\} \ .
\end{eqnarray*}
\end{cor}
We apply the Corollary to minimal surfaces $M^2 \subset X^3(\kappa)$ in
a $3$-dimensional space of constant curvature $\kappa$. 
The second fundamental form is a
Codazzi tensor. The Gauss equation  $ S = 2 \kappa - 2 a^2$ yields finally 
the result
\begin{eqnarray*}
\lambda_1^2  \,   \geq   \,   \kappa \, + \,  \inf_M \Big(
(\overline{\lambda}_1^2  \, -\, 1) \, a^2\, + \, 
a^4 \, \Delta(a^{-4}) \Big) \ .
\end{eqnarray*}
%
\section{Deformation of the metric via a Codazzi tensor}
\noindent
In this section we establish some lemmata that we will need later to 
prove Theorems $1.1$ and $1.2$. Consider a nondegenerate symmetric 
$(0, 2)$-tensor field $\beta$ on $(M^n, g)$ and define a
new metric $\overline{g}$ by 
\begin{equation}  
\overline{g}(X, Y) = g ( \beta(X), \beta(Y) ).      
\end{equation}
The Levi-Civita connection $\overline{\nabla}$ of $(M^n,
\overline{g})$ is related to the Levi-Civita connection $\nabla$ of 
$(M^n, g)$  by \cite{Kim2002}
\begin{equation}
\overline{\nabla}_{\beta^{-1}(X)} \left( \beta^{-1}(Y) \right) \ =
\ \beta^{-1} \left( \nabla_{\beta^{-1}(X)} Y \right) +  \beta^{-1}
\left( \Lambda (X, Y) \right) ,
\end{equation}
where $ \Lambda$ is the $(1,2)$-tensor field defined by
\begin{eqnarray}
2 \, g ( \Lambda (X , Y ),  Z) &  =  &   g \left(  Z  , \, \beta
\{  ( \nabla_{\beta^{-1}(X)} \beta^{-1} ) (Y) \}  -
\beta \{  ( \nabla_{\beta^{-1}(Y)}  \beta^{-1} )(X) \}  \right)   \nonumber   \\
&     &   \nonumber  \\
&    & +  g \left(  Y  , \, \beta \{  ( \nabla_{\beta^{-1}(Z)}
\beta^{-1} )(X) \}  -
\beta \{  ( \nabla_{\beta^{-1}(X)}  \beta^{-1} )(Z) \}  \right)      \nonumber   \\
&     &   \nonumber  \\
&     & +  g \left(  X  , \, \beta \{  ( \nabla_{\beta^{-1}(Z)}
\beta^{-1} )(Y) \}  - \beta \{  ( \nabla_{\beta^{-1}(Y)}  \beta^{-1} )(Z)
\} \right).
\end{eqnarray}
Note that the tensor $\Lambda$ satisfies
\[ g( \Lambda(X, Z), Y) + g( \Lambda(X, Y), Z) = 0  \]
 for all vector fields $X, Y, Z$.
Using  formula (12), we can relate the Riemann curvature tensor
$\overline{R}$ of $(M^n, \overline{g})$ to the one $R$
of $(M^n, g)$ by
\begin{eqnarray*}
&     & \overline{R} (\beta^{-1}X , \beta^{-1}Z)(\beta^{-1}Y) - 
\beta^{-1} \{ R( \beta^{-1}X, \beta^{-1}Z)(Y) \}    \\
&      &      \\
&    =  &  \beta^{-1} \{ (\nabla_{\beta^{-1}(X)} \Lambda)(Z, Y) - 
( \nabla_{\beta^{-1}(Z)} \Lambda)(X, Y) \}   +  \beta^{-1} 
\{ \Lambda(X, \Lambda(Z, Y))     \nonumber     \\
&      &        \\
&      & - \Lambda(Z, \Lambda(X, Y)) \} + \beta^{-1}  
\{ \Lambda ( \Lambda (Z, X) - \Lambda(X, Z), \, Y)  \}.
\end{eqnarray*}
Let $( E_1 , \ldots , E_n )$ be a local $g$-orthonormal frame
field on $( M^n , g )$.
Then the scalar curvature $\overline{S}$ of $(M^n, \overline{g})$ is 
expressed as
\begin{eqnarray}
&     & \overline{S} - \sum_{i,j=1}^n g( E_i , 
R( \beta^{-1}E_i, \beta^{-1}E_j)(E_j) )    \nonumber    \\
&      &    \nonumber   \\
&   =   &  2 \sum_{i,j=1}^n g( E_i, \, (\nabla_{\beta^{-1}(E_i)} 
\Lambda)(E_j, E_j)) - \sum_{i,j,k=1}^n \Lambda_{iik} \Lambda_{jjk}   
 - \sum_{i,j,k=1}^n \Lambda_{ijk}  \Lambda_{jik},
\end{eqnarray}
where $\Lambda_{ijk} : = g ( \Lambda(E_i, E_j) , E_k )$. We now
review briefly the behavior of the Dirac operator under the
deformation (11) of metrics.  Let $\Sigma (M)_g$ and $\Sigma
(M)_{\overline{g}}$ be the spinor bundles of $(M^n, g)$ and $(M^n,
\overline{g})$, respectively. There are natural isomorphisms
$\beta^{-1} : T(M)  \longrightarrow T(M)$ and
${\widehat{\beta^{-1}}} : \Sigma {(M)}_g \longrightarrow \Sigma
{(M)}_{\overline{g}} $   preserving the inner products of vectors
and spinors as well as the Clifford multiplication:
\begin{eqnarray*}
&     &  \overline{g} ( \beta^{-1} X, \beta^{-1} Y )  = g (X , Y) ,  
\qquad  \langle \widehat{\beta^{-1}} \varphi,
\widehat{\beta^{-1}} \psi \rangle_{\, \overline{g}} \, = \, \langle \varphi, 
\psi \rangle_g,      \\
&       &       \\
&     &    (\beta^{-1} X)  \cdot ( \widehat{\beta^{-1}} \psi ) = 
\widehat{\beta^{-1}} ( X \cdot \psi ) ,   \qquad
X, Y  \in \Gamma( T(M)) , \quad \varphi , \psi \in \Gamma( \Sigma (M)_g ) .
\end{eqnarray*}
For each spinor field $\psi$ on $(M^n, g)$ we denote by
$\overline{\psi} := \widehat{\beta^{-1}} ( \psi )$ the
corresponding spinor field on $(M^n, \overline{g})$. We will use the
same notation for vector fields, $\overline{X} := \beta^{-1}(X)$.
It follows from (12) that the spinor  derivatives
$\overline{\nabla}, \, \nabla$ are related by
\begin{equation}
  \overline{\nabla}_{\beta^{-1}(E_j)}   \overline{\psi}
  =    \overline{\nabla_{\beta^{-1}(E_j)} \psi}  + \frac{1}{4} \sum_{k,l=1}^n  
\Lambda_{jkl} \overline{E_k} \cdot \overline{E_l} \cdot \overline{\psi}  .
\end{equation}
\noindent Let $\omega$ and $\Omega$ be a 1-form and a 3-form generated by the 
tensor $\Lambda$ via
\[ \omega= \sum_{j,k=1}^n \Lambda_{jjk} E^k,\qquad E^k := g (\cdot,E_k),\]
and
\[ \Omega = \sum_{j<k<l} (\Lambda_{jkl} + \Lambda_{klj} + \Lambda_{ljk}) E^j 
\wedge E^k \wedge E^l , \]
respectively.  The Dirac operator $\overline{D}$ of $(M^n, \overline{g})$ 
can be  expressed through the $\beta$-twist $D_\beta$ of $D$ as
\begin{eqnarray}
\overline{D} \, \overline{\psi} & = & \sum_{i=1}^n \overline{E_i} \cdot 
\overline{\nabla}_{\overline{E_i}} \, \overline{\psi}      
\ =  \  \sum_{i=1}^n   \overline{E_i \cdot \nabla_{\beta^{-1}(E_i)} \psi} 
+ \frac{1}{4}
\sum_{j,k,l=1}^n   \Lambda_{jkl} \overline{E_j \cdot E_k \cdot E_l \cdot \psi} 
\nonumber  \\
& = &  \overline{D_{\beta} \psi}
- \frac{1}{2} \, \overline{\omega \cdot \psi}  + \frac{1}{2} \, 
\overline{\Omega \cdot \psi}
\end{eqnarray}
and the square $\overline{D}^2$ of the Dirac operator $\overline{D}$ by
\begin{eqnarray}
\overline{D}^2 \, \overline{\psi}   & = &  \overline{( D_{\beta} \circ 
D_{\beta}) (\psi)}
- \frac{1}{2} \, \overline{\omega \cdot D_{\beta} \psi}  + \frac{1}{2} \, 
\overline{\Omega \cdot D_{\beta} \psi}  \nonumber  \\
&   & - \frac{1}{2} \, \overline{D_{\beta} (\omega \cdot \psi)} 
+ \frac{1}{2} \, \overline{D_{\beta} (\Omega \cdot \psi)}
+ \frac{1}{4} \, \overline{\omega \cdot \omega \cdot \psi} \nonumber \\
&     &  - \frac{1}{4} \, \overline{\Omega \cdot \omega \cdot \psi} 
- \frac{1}{4} \, \overline{\omega \cdot \Omega \cdot
\psi} + \frac{1}{4} \, \overline{\Omega \cdot \Omega \cdot \psi}.
\end{eqnarray}
\noindent
In this paper we focus our attention on an interesting property of the 
tensor $\Lambda$. Since
\begin{eqnarray*}
  \Lambda (X, Y) - \Lambda (Y, X) & = &  \beta \{ (\nabla_{\beta^{-1}X} 
\beta^{-1})(Y) \}-\beta \{ (\nabla_{\beta^{-1}Y} \beta^{-1})(X) \} \\
& =   &  -  (\nabla_{\beta^{-1}X} \beta) (\beta^{-1}Y) + (\nabla_{\beta^{-1}Y}
\beta) (\beta^{-1}X),
\end{eqnarray*}
we observe that $\Lambda \equiv 0$. Consequently, all the equations 
simplify remarkably when $\beta$ is a Codazzi tensor.
\begin{lem}
Let $\beta$ be a nondegenerate Codazzi tensor on $(M^n, g)$. Then
we have:
\begin{eqnarray}
\overline{S} & = & \sum_{i,j=1}^n  g ( E_i,  R( \beta^{-1}E_i, \beta^{-1}E_j)
(E_j) ), \\
\overline{\nabla}_X \overline{\psi}  &  =  &  \overline{\nabla_X \psi}, \\
\overline{D} \, \overline{\psi} &  = & \overline{D_{\beta} \psi},\\
\overline{D}^2 \, \overline{\psi} &  = & \overline{(D_{\beta} \circ D_{\beta})
( \psi)} .
\end{eqnarray}
\end{lem}
We close the section with  some more lemmata needed in
the next section.
\begin{lem}
Let $\beta$ be a nondegenerate symmetric tensor field on $(M^n, g
)$. If there exists a nontrivial spinor field $\psi$ on $(M^n, g)$
such that
\begin{equation}
\nabla_X \psi = p \, X \cdot D \psi + q \, \beta^{-1}(X) \cdot
D_{\beta} \psi
\end{equation}
holds for some real-valued functions $p , q : M^n  \longrightarrow \mathbb{R}$ 
and for all vector fields $X$, then
\begin{eqnarray}
&    &  (1 + np) D \psi =  - q  \, {\rm tr}(\beta^{-1}) D_{\beta}
\psi,  \\
&     &   \nonumber     \\
 &    &  (1 + q \vert \beta^{-1} \vert^2
) D_{\beta} \psi = - p \, {\rm tr}(\beta^{-1}) D \psi.
\end{eqnarray}
\end{lem}
\begin{lem}
Let $( ,  ) := {\rm Re} \langle ,  \rangle$ denote the real part of the 
standard Hermitian  product $\langle , \rangle$ on the spinor bundle 
$\Sigma (M)$ over $M^n$. Let $\psi$ and $F$ be a spinor field and a 
real-valued function on $M^n$, respectively.
Then we have
\[  F \cdot  \triangle ( \psi , \psi ) - ( \psi , \psi ) \cdot 
\triangle F = {\rm div} \{  ( \psi , \psi ) {\rm grad} F  -  F  \, {\rm grad} 
( \psi , \psi ) \} .   \]
\end{lem}
%
\section{Proof of the Theorems}
%
Note that the volume form $\overline{\mu}$ of $( M^n , \overline{g} )$ is 
related to the one $\mu$ of $( M^n , g )$ by
\begin{equation}
\overline{\mu} = \vert {\rm det}( \beta^{-1} )  \vert  \mu .
\end{equation}
Let $\mathcal{Q} : \Gamma ( T(M) ) \times \Gamma( \Sigma (M)_g
 ) \longrightarrow \Gamma ( \Sigma(M)_g )$ be a {\it twistor-like
operator} defined by
\[  \mathcal{Q}_X (\varphi) = \nabla_X \varphi - p \, X \cdot D
\varphi - q \, \beta^{-1}(X) \cdot D_{\beta} \varphi ,  \]
where $p, q : M^n  \longrightarrow \mathbb{R}$ are some real-valued functions.
Then we have
\begin{eqnarray}
\sum_{j=1}^n ( \mathcal{Q}_{E_j} (\varphi)\!\!\!\! &,& \!\!\!\!
\mathcal{Q}_{E_j} (\varphi)) \ =\ \nonumber \\
&&\!\!\!\!
 {\rm div} \big[ \sum_{j=1}^n ( \varphi , \, E_j  \cdot D \varphi + \nabla_{E_j} \varphi ) E_j    \big]
 + ( n p^2 + 2p + 1 ) ( D \varphi,  D \varphi ) \nonumber \\
&&\!\!\!\!
 - \frac{1}{4} S ( \varphi , \varphi) + 
\{ q^2 \vert \beta^{-1} \vert^2 + 2q \} ( D_{\beta} \varphi , 
D_{\beta} \varphi )
 + 2pq \, {\rm tr}( \beta^{-1}) ( D \varphi, D_{\beta} \varphi ).
\end{eqnarray}
Now, let $\varphi = \psi$ be an eigenspinor of $D$ with eigenvalue $\lambda
\not= 0 \in \mathbb{R}$. By Lemma 2.3, we then see that
\begin{equation}
\int\limits_{M^n}  F {\rm div} \bigg[ \sum_{j=1}^n ( \psi , \, E_j  \cdot D 
\psi + \nabla_{E_j} \psi ) E_j    \bigg] \mu
=  - \frac{1}{2} \int_{M^n}  ( \psi, \psi ) \triangle (F) \mu
\end{equation}
holds for any real-valued function 
$F : M^n \longrightarrow \mathbb{R}$, since $( \psi , E_j \cdot \psi ) = 0$ 
vanishes identically.
Let $\overline{\lambda}_1 \not= 0$ be the smallest eigenvalue of 
$\overline{D}$.
Making use of (20), (25), (26), (27) and introducing free functions 
$F, B, C : M^n \longrightarrow \mathbb{R}$ (we assume that $ F$ and $B^2$ 
are positive functions.) to control
the unnecessary terms, we  compute
\begin{eqnarray}
H_1   &  :=  &  \int\limits_{M^n} \left[  ( \overline{D} \, \overline{\psi} , \,  
\overline{D} \, \overline{\psi} ) -  \overline{\lambda}_1^2 ( \overline{\psi}, 
\overline{\psi} )  \right] \    \overline{\mu}   \nonumber \\
& &  +  \int\limits_{M^n}  \left[   F
\sum_{j=1}^n ( \mathcal{Q}_{E_j} (\psi) ,  \mathcal{Q}_{E_j}(\psi) ) 
+  B^2 (  D_{\beta} \psi - C D \psi,   D_{\beta} \psi - C D \psi )  \right] 
\mu    \nonumber  \\
&  =  &   \int\limits_{M^n}  \left[   \lambda^2 \Big((n p^2 + 2p +1) F + B^2
  C^2  \Big)  - \frac{1}{4} F S  -  \overline{\lambda}_1^2  \,  
\vert \det( \beta^{-1} )  \vert  - \frac{1}{2} (\triangle F) \right]
(\psi, \psi ) \, \mu  \nonumber  \\
&  & + \int\limits_{M^n}  \bigg[2 \lambda \Big(  pq F \, {\rm tr}(\beta^{-1}) -  B^2 C \Big) ( \psi, D_{\beta} 
\psi ) \nonumber\\
& & + \, \Big( (q^2 \vert \beta^{-1} \vert^2 + 2q)F + B^2  +  \vert {\rm det}( 
\beta^{-1} )  \vert \Big) ( D_{\beta} \psi, D_{\beta} \psi )        
\bigg]  \mu .
\end{eqnarray}
We choose the functions $B, C$ in such a way that the 
second integral of (28) vanishes
and the equations (23), (24) are satisfied with $D_{\beta} \psi
= C D \psi$.  To this end, it is required that the relations
\begin{equation}
B^2 = - q F ( 1 + q \vert \beta^{-1} \vert^2 ) \ , \quad
C^2 = \frac{p (1+ np)}{q ( 1 + q \vert \beta^{-1} \vert^2 )}
\end{equation} 
hold.
Note that (29) implies, in particular, the restriction (4) :
\[  - \frac{1}{\vert \beta^{-1} \vert^2 } < q < 0 ,  \qquad  - \frac{1}{n} < p < 0   .    \]
Now choose
\begin{equation}   F = - \
\frac{ \vert {\rm det}( \beta^{-1} )  \vert }{ q} ,  \qquad  B^2 =   
\vert {\rm det}( \beta^{-1} )  \vert
( 1 + q \vert \beta^{-1} \vert^2 )    \end{equation}
so that the last line in the latter part of  (28) vanishes.
Then we obtain
\begin{eqnarray*}
\int\limits_{M^n} \left[\lambda^2(p+1) F  - \frac{1}{4} 
F S  - \overline{\lambda}_1^2 \, \vert \det( \beta^{-1} )
  \vert   - \frac{1}{2} (\triangle F) 
\right] (\psi, \psi ) \,  \mu   \geq  0 \, ,
\end{eqnarray*}
which proves the inequality of Theorem $1.1$.   The functions
 optimal for  $p$ and $q$ are to be found when considering
the limiting case.  
The former part of (28) yields in the limiting case that
\begin{equation}
\overline{D} \, \overline{\psi_1} = \overline{\lambda}_1 \overline{\psi_1} 
= \overline{D_{\beta} \psi_1} =  C \lambda_1  \, \overline{\psi_1}. 
\end{equation}
Since $\overline{\lambda}_1  = C \lambda_1$, we find that the function 
$C$  must be a nonzero constant
$C= c \not= 0  \in {\mathbb R}$.  Then, from (31) and (23), (24), we obtain
the two relations in (5) immediately. The condition (6), (7) 
for the limiting case is easy to check.

To prove Theorem $1.2$ we consider the integral
\begin{eqnarray}
H_2   &  :=  &  \int\limits_{M^n}\left[(\overline{D} \, \overline{\psi}, \,  
\overline{D} \, \overline{\psi} ) -  \overline{\lambda}_1^2 ( \overline{\psi}, 
\overline{\psi} )  \right] \    \overline{\mu}      
+  \int\limits_{M^n}  \bigg[ F\sum_{j=1}^n ( \mathcal{Q}_{E_j} (\psi) ,  
\mathcal{Q}_{E_j}(\psi) )  \bigg] \mu \nonumber \\
&  =  &   \int\limits_{M^n}  \left[\lambda^2 (n p^2 + 2p +1) F  
- \frac{1}{4} F S   
- \overline{\lambda}_1^2  \,  \vert {\rm det}( \beta^{-1} )  \vert   
- \frac{1}{2} (\triangle F)  \right]  (\psi, \psi ) \, \mu \nonumber \\
&     &  
+   \int\limits_{M^n} \left[ (q^2 \vert \beta^{-1} \vert^2 + 2q)F  
+  \vert {\rm det}( \beta^{-1} )  \vert \right] ( D_{\beta} \psi, D_{\beta} 
\psi )  \, \mu 
\end{eqnarray}
and choose the free parameters $p, q, F$ as
\[ 
 p = - \frac{1}{n} , \qquad  q =  - \frac{1}{ \vert  \beta^{-1}  
\vert^2 } ,  \qquad F =  \vert {\rm det}(
\beta^{-1} )  \vert  \vert \beta^{-1} \vert^2  . \]
Then the third line in the latter part of (32) vanishes and we have
\[
\int\limits_{M^n}  \left[\frac{(n-1) \lambda^2 F}{n}     -
  \frac{1}{4} F S  - \overline{\lambda}_1^2    
\vert \det( \beta^{-1} )  \vert   - 
\frac{1}{2} (\triangle F)  \right]  (\psi, \psi ) \, \mu \geq 0 \, . \]
This proves the inequality (8).  The condition (10) for the 
limiting case is clear.
\begin{NB} Let $\overline{\lambda}_1 \not= 0  \in {\mathbb R}$ be the smallest
 eigenvalue of $\overline{D}$.  Suppose that there exist a nonzero constant
$\lambda \not= 0 \in {\mathbb R}$ and a spinor field $\psi$ such that  
the following equations hold:
\begin{eqnarray*}
D \psi \ = \ \lambda \psi \ , \quad D_{\beta} \psi  \ =  \
\overline{\lambda}_1 \psi  \ ,   \quad \nabla_X \psi =   \lambda \, p \, X 
\cdot \psi  +
\overline{\lambda}_1 q \, \beta^{-1}(X) \cdot \psi .   
\end{eqnarray*}
Then it turns out that $\lambda = \lambda_1$ is equal to the smallest 
eigenvalue of $D$ and, in the limiting case, the constant $c$ in  (5) is
related to $\lambda_1, \, \overline{\lambda}_1$ by 
$\overline{\lambda}_1 = c \, \lambda_1$.
\end{NB}

\end{document}